\documentclass[11pt,
              ]{amsart}

\usepackage{amssymb}
\usepackage{amsxtra}
\usepackage[mathscr]{eucal}
\usepackage{graphics}

\numberwithin{equation}{section}

\newtheorem{theorem}{Theorem}[section]

\newtheorem{lemma}[theorem]{Lemma}
\newtheorem{corollary}[theorem]{Corollary}
\newtheorem{conjecture}[theorem]{Conjecture}

\theoremstyle{definition}

\newtheorem{definition}[theorem]{Definition}
\newtheorem{example}[theorem]{Example}

\newtheorem{remark}[theorem]{Remark}

\newcommand{\ghat}{\widehat{G}}
\newcommand{\lam}{\lambda}

\newcommand{\RR}{\mathbb R}
\newcommand{\m}{\mathsf{m}}
\newcommand{\tor}{\mathsf{t}}
\newcommand{\TT}{\mathbb{T}}
\newcommand{\ZZ}{\mathbb{Z}}
\newcommand{\z}[1]{\ZZ_{\!/#1}}

\renewcommand{\ge}{\geqslant}
\renewcommand{\le}{\leqslant}

\begin{document}

\title[Lehmer's Problem]
  {Lehmer's Problem for\\ Compact Abelian Groups}

\author[Douglas Lind]{Douglas Lind}

\address{Department of Mathematics,
   Box 354350, University of Washington, Seattle, WA 98195--4350}

\email{lind@math.washington.edu}

\date{\today}




\begin{abstract}
  We formulate Lehmer's Problem about the Mahler measure of
  polynomials for general compact abelian groups, introducing a
  Lehmer constant for each such group. We show that all
  nontrivial connected compact groups have the same Lehmer
  constant, and conjecture the value of the Lehmer constant for
  finite cyclic groups. We also show that if a group has
  infinitely many connected components then its Lehmer constant
  vanishes. 
\end{abstract}

\maketitle

\section{Introduction}\label{sec:introduction}
  
Let $f\in\ZZ[x^{\pm1}]$ be a Laurent polynomial with integer
coefficients. Define its \textit{logarithmic Mahler measure} to
be
\begin{displaymath}
   \m(f)=\int_0^1\log|f(e^{2\pi i s})|\,ds.
\end{displaymath}
Lehmer's Problem asks whether $\m(f)$ can be arbitrarily small
but positive. Equivalently, does
\begin{displaymath}
  \inf\{\, \m(f):f\in\ZZ[x^{\pm1}],\ \m(f)>0 \,\} 
\end{displaymath}
equal zero? The smallest positive value of $\m(f)$ known, found
by Lehmer himself \cite{Lehmer}, is attained by
\begin{equation}\label{eqn:lehmer-polynomial}
   f_L(x)=x^{10}+x^9-x^7-x^6-x^5-x^4-x^3+x+1,
\end{equation}
for which $\m(f_L)\approx 0.16235$. For accounts of this problem,
see \cite{Boyd1}, \cite{Boyd2}, and for connections with dynamics
see \cite{EW}. 

Lehmer's Problem can be formulated for arbitrary compact abelian
groups. To do so, let $G$ be a compact abelian group with
normalized Haar measure~ $\mu$. Let $\ghat$ denote its
(multiplicative) dual group of characters, and $\ZZ[\ghat]$ be
the ring of integral combinations of characters. For
$f\in\ZZ[\ghat]$ define its \textit{logarithmic Mahler measure
over $G$} to be
\begin{displaymath}
   \m(f)=\m_G(f)=\int_G\log|f|\,d\mu.
\end{displaymath}
Since $\log 0 =-\infty$, if $f$ vanishes on a set of positive
$\mu$-measure then $\m(f)=-\infty$, but otherwise $\m(f)\ge0$ by
Lemma \ref{lem:positivity} below.

\begin{definition}\label{def:lehmer-constant}
   The \textit{Lehmer constant} of a compact abelian group $G$ is
   \begin{displaymath}
      \lam(G)=\inf\{\,\m_G(f):f\in\ZZ[\ghat],\ \m_G(f)>0\,\}.
   \end{displaymath}
\end{definition}

Lehmer's Problem therefore asks whether $\lam(\TT)=0$, where
$\TT=\RR/\ZZ$. 

We show in Theorem \ref{thm:connected} that all nontrivial
connected groups have the same Lehmer constant $\lam(\TT)$. On
the other hand, Theorem \ref{thm:disconnected} shows that if $G$
has infinitely many connected components, then $\lam(G)=0$. We
compute the Lehmer constant of some finite groups, and conjecture
a simple formula for cyclic groups and also for products of
two-element groups. Finally, we show that if $\lam(\TT)=0$, then
the only groups which have positive Lehmer constant are finite.

\section{Preliminaries}\label{sec:preliminaries}

Let $G$ be a compact abelian group, and let $f\in\ZZ[\ghat]$. If
$f$ vanishes on a set of positive $\mu$-measure then clearly
$\m(f)=-\infty$. 

\begin{lemma}\label{lem:positivity}
   Let $f\in\ZZ[\ghat]$. If $\mu(\{x:f(x)=0\})=0$, then
   $\m(f)\ge0$. 
\end{lemma}

\begin{proof}
   The characters appearing in $f$ generate a subgroup $\Delta$
   of $\ghat$. Let $H=\Delta^\perp$ and $\pi\colon G\to G/H$ be
   the quotient map. Then $f$ is constant on cosets of $H$, and
   so defines $\widetilde{f}\in\ZZ[(G/H)\sphat\,]=\ZZ[\Delta]$ such
   that $\widetilde{f}\circ\pi=f$. Clearly
   $\m_{G/H}(\widetilde{f})=\m_G(f)$, and $\widetilde{f}$
   vanishes on a set of positive measure if and only if $f$
   does. Hence we may assume that $\ghat$ is finitely generated.

   Let $\z{n}$ denote $\ZZ/n\ZZ$.
   Applying the structure theorem for finitely generated abelian
   groups to $\ghat$, we see that $G$ is isomorphic to
   \begin{displaymath}
      \z{n_1}\oplus\dots\oplus\z{n_r}\oplus\TT^k
   \end{displaymath}
   for suitable integers $n_1,\dots,n_r$ and $k$. Consider
   $f\in\ZZ[\ghat]$ as a function of $r+1$ variables
   $j_1\in\z{n_1}$, $\dots$, $j_r\in\z{n_r}$, and
   $\mathbf{s}\in\TT^k$. Since the set where $f$ vanishes is
   null, each function
   \begin{displaymath}
      f_{j_1\dots j_r}(\mathbf{s})=f(j_1,\dots,j_r,\mathbf{s})
   \end{displaymath}
   is a nonvanishing complex combination of characters on $\TT^k$,
   and hence by \cite[Lemma 3.7]{EW} we see that
   $\m_{\TT^k}(f_{j_1\dots j_r})>-\infty$.
   Thus
   \begin{align*}
      \m(f) &= \frac{1}{n_1\dots n_r} \sum_{j_1=0}^{n_1-1}\dots
      \sum_{j_r=0}^{n_r-1} \int_{\TT^k}\log|f_{j_1\dots j_r}
      (\mathbf{s})|\,d\mathbf{s} \\
      &= \frac{1}{n_1\dots n_r} \m_{\TT^k}(g) >-\infty,
   \end{align*}
   where
   \begin{displaymath}
      g(\mathbf{s})=\prod_{j_1=0}^{n_1-1}\dots\prod_{j_r=0}^{n_r-1}
      f_{j_1\dots j_r}(\mathbf{s}).
   \end{displaymath}
   
   Since for each $i$ the product over $j_i$ covers all the
   $n_i$th roots of unity, the coefficients of $g$ are algebraic
   integers fixed by all elements in the Galois group of the
   field they generate, and so $g\in\ZZ[(\TT^k)\sphat\,]$. Hence
   $\m_{\TT^k}(g)\ge0$ by \cite[Lemma 3.7]{EW}, showing that
   $\m(f)\ge0$.
\end{proof}

We repeatedly use the following observation.

\begin{lemma}\label{lem:quotient}
   Let $H$ be a closed subgroup of $G$. Then
   $\lam(G)\le\lam(G/H)$. 
\end{lemma}

\begin{proof}
   Let $\pi\colon G\to G/H$ be the quotient map. If
   $f\in\ZZ[(G/H)\sphat\,]$, then $f\circ\pi\in\ZZ[\ghat]$ and
   $\m_G(f\circ\pi)=\m_{G/H}(f)$. The conclusion then follows from
   the definition of $\lam$. 
\end{proof}

\section{Connected groups}\label{sec:connected}

In this section we prove that all nontrivial connected groups
have the same Lehmer constant.

\begin{lemma}\label{lem:T^n}
   $\lam(\TT^k)=\lam(\TT)$ for all $k\ge1$.
\end{lemma}

\begin{proof}
   Since $\TT\cong\TT^k/\TT^{k-1}$, Lemma \ref{lem:quotient}
   shows that $\lam(\TT^k)\le\lam(\TT)$.

   To prove the reverse inequality, fix $\epsilon>0$. Choose
   $f\in\ZZ[(\TT^k)\sphat\,]$ with
   \begin{displaymath}
      0<\m_{\TT^k}(f)<\lam(\TT^k)+\epsilon.
   \end{displaymath}
   For $\mathbf{r}=(r_1,\dots,r_k)\in\ZZ^k$ define
   $f_{\mathbf{r}}(s)=f(r_1s,\dots,r_ks)$, so that
   $f_{\mathbf{r}}\in\ZZ[\widehat{\TT}]$. By \cite{Lawton} we can
   find $\mathbf{r}$ so that
   \begin{displaymath}
      |\m^{}_{\TT}(f_{\mathbf{r}})-\m_{\TT^k}(f)|
      <\min\{\m_{\TT^k}(f),\epsilon\}.
   \end{displaymath}
   It follows that $\lam(\TT)\le\lam(\TT^k)+2\epsilon$. Since
   $\epsilon>0$ was arbitrary, we obtain that
   $\lam(\TT)\le\lam(\TT^k)$.
\end{proof}

\begin{theorem}\label{thm:connected}
   If $G$ is a nontrivial connected compact abelian group, then
   $\lam(G)=\lam(\TT)$. 
\end{theorem}

\begin{proof}
   Since $G$ is nontrivial and connected, it has a quotient
   isomorphic to ~$\TT$, and hence $\lam(G)\le\lam(\TT)$ by Lemma
   \ref{lem:quotient}.

   To prove the reverse inequality, observe that $\ghat$ is
   torsion-free. Let $f\in\ZZ[\ghat]$ with $\m_G(f)>0$. The
   subgroup $\Delta$ of $\ghat$ generated by the characters in
   $f$ is therefore isomorphic to $\ZZ^k$ for some $k\ge1$, and
   so $\widehat{\Delta}=G/\Delta^\perp\cong\TT^k$. As in the
   proof of Lemma \ref{lem:positivity}, $f$ induces
   $\widetilde{f}\in\ZZ[\Delta]$ with
   $\m_{\widehat{\Delta}}(\widetilde{f})=\m^{}_G(f)$. By Lemma
   \ref{lem:T^n}, we have that
   \begin{displaymath}
      \m^{}_G(f)=\m_{\widehat{\Delta}}(\widetilde{f})\ge\lam(\TT^k)
      =\lam(\TT),
   \end{displaymath}
   proving that $\lam(G)\ge\lam(\TT)$.
\end{proof}

\section{Finite groups}\label{sec:finite}

In this section we consider the problem of computing the Lehmer
constant of some finite abelian groups.

Trivially $\lam(\z{1})=\log 2$.

Define characters $\chi_k$ on $\z{n}$ by $\chi_k(j)=e^{2\pi i j
k/n}$. Thus the character group is
$\widehat{\z{n}}=\{\chi_0,\chi_1,\dots,\chi_{n-1}\}$.

\begin{example}\label{exam:Z/2}
   Let $G=\z{2}$. For $f=a\chi_0+b\chi_1\in\ZZ[\widehat{\z{2}}] $ we
   have that
   \begin{displaymath}
      \m(f)=\frac12 \bigl(\log|a+b|+\log|a-b|\bigr)
      =\frac12 \log|a^2-b^2|.
   \end{displaymath}
   It is easy to check that the smallest value of $|a^2-b^2|$
   greater than 1 is 3, so that
   $\displaystyle\lam(\z{2})=\frac12\log3$. 
\end{example}

\begin{example}\label{exam:Z/3}
   Let $G=\z{3}$. A simple calculation shows that
   \begin{displaymath}
      \m(a\chi_0+b\chi_1+c\chi_2)=\frac13 \log|a^3+b^3+c^3-3abc|.
   \end{displaymath}
   Choosing $a=b=1$ and $c=0$ gives a value of 2 for the
   expression inside the absolute value, which is clearly the
   smallest possible value greater than~ 1. Hence
   $\displaystyle\lam(\z{3})=\frac13\log2$.
\end{example}

\begin{example}\label{exam:Z/4}
   Let $G=\z{4}$. Then
   \begin{displaymath}
      \m(a\chi_0+b\chi_1+c\chi_2+d\chi_3)=\frac14
      \log\bigl| (a+b+c+d)(a-b+c-d)[(a-c)^2+(b-d)^2]\bigr|.
   \end{displaymath}
   Putting $a=b=c=1$ and $d=0$ gives $\displaystyle\frac14\log3$.

   Suppose the product inside the absolute value were $\pm2$. Then
   the factorization would be either $2\cdot1\cdot1$,
   $1\cdot2\cdot1$, or $1\cdot1\cdot2$. The first two
   factorizations cannot occur since the difference of the first
   two factors must be even. Therefore $|a+c|\le2$ and
   $(a-c)^2+(b-d)^2=2$. A straightforward search rules out all
   possibilities to attain $\pm2$. Hence
   $\displaystyle\lam(\z{4})=\frac14\log3$. 
\end{example}

There is a simple upper bound for the Lehmer constant of a finite
group. 

\begin{lemma}\label{lem:upper-bound}
   Let $F$ be a finite abelian group with cardinality
   $|F|\ge3$. Then
   \begin{equation}\label{eqn:upper-bound}
      \lam(F)\le\frac1{|F|} \log\bigl(|F|-1\bigr).
   \end{equation}
\end{lemma}

\begin{proof}
   Let $\chi_0\in\widehat{F}$ be the trivial character, and
   $\delta_0$ be the unit mass at $0\in F$. Then
   \begin{displaymath}
      g=\sum_{\chi\in\widehat{F}} \chi=|F|\delta_0.
   \end{displaymath}
   Put $f=g-\chi_0$, so that
   \begin{displaymath}
      f(x)=
      \begin{cases}
         |F|-1 &\text{if $x=0$},\\
         -1    &\text{if $x\ne0$}.
      \end{cases}
   \end{displaymath}
   Then $\displaystyle
   \m_F(f)=\frac1{|F|}\log\bigl(|F|-1\bigr)>0$ since $|F|\ge3$,
   proving \eqref{eqn:upper-bound}.
\end{proof}

To estimate the Lehmer constant of a cyclic group, it is
convenient to introduce the following arithmetical function.

\begin{definition}\label{def:prime}
   For an integer $n\ge2$ let $\rho(n)$ denote the smallest prime
   number that does not divide $n$.
\end{definition}

It is easy to see by using the fact that $\prod_{p \le x} p =
\exp(x(1 + o(x))$ (equivalent to the Prime Number Theorem) that
$\rho(n) \le \log n (1 + o(1))$ as $n \to \infty$.

\begin{theorem}\label{thm:finite-bound}
   For all integers $n\ge2$ we have that
   \begin{displaymath}
      \lam(\z{n})\le\frac1n \log\rho(n).
   \end{displaymath}
\end{theorem}

\begin{proof}
   Let $p=\rho(n)$, and $\Phi_p$ be the $p$th cyclotomic
   polynomial. Put $f_p(j)=\Phi_p(e^{2\pi i j/n})$, so that
   $f_p\in\ZZ[\widehat{\z{n}}]$. Then
   \begin{displaymath}
      \m(f_p)=\frac1n\log|R\bigl(\Phi_p(x),x^n-1\bigr)|,
   \end{displaymath}
   where $R(\cdot,\cdot)$ denotes the resultant. Since $p$ does
   not divide $n$, if follows from \cite{Apostol} that
   $R\bigl(\Phi_p(x),x^n-1\bigr)=p=\rho(n)$.
\end{proof}

\begin{corollary}\label{cor:n-odd}
   If $n$ is odd then $\displaystyle\lam(\z{n})=\frac1n\log2$.
\end{corollary}

\begin{proof}
   Since $\rho(n)=2$ we see that
   $\displaystyle\lam(\z{n})\le\frac1n\log2$ by Theorem
   \ref{thm:finite-bound}. An argument similar to that in the
   proof of Lemma \ref{lem:positivity} shows the reverse inequality.
\end{proof}

Note that the upper bound in Theorem \ref{thm:finite-bound} is
actually the correct value of $\lam(\z{n})$ for $n=2,3$, and $4$
computed in Examples \ref{exam:Z/2}, \ref{exam:Z/3}, and
\ref{exam:Z/4} above. Further computational evidence suggests the
following.

\begin{conjecture}\label{conj:cyclic}
   $\displaystyle\lam(\z{n})=\frac1n\log\rho(n)$ for all $n\ge2$.
\end{conjecture}

Next, consider groups of the form
$\z{2}\oplus\dots\oplus\z{2}=\z{2}^n$.

\begin{example}\label{exam:Z/2^2}
   Let $G=\z{2}\oplus\z{2}$. The characters
   $\chi_{ij}(x,y)=\chi_i(x)\chi_j(y)$ form the dual group
   $\ghat=\{\chi_{00},\chi_{01},\chi_{10},\chi_{11}\}$. Then
   \begin{displaymath}
      \m\Bigl( \sum_{i,j=0}^1 a_{ij}\chi_{ij} \Bigr) =
      \frac14 \log \Bigl| \prod_{r,s=0}^1 \Bigl(
      \sum_{i,j=0}^1 (-1)^{ir+js}a_{ij}\Bigr)\Bigr|.
   \end{displaymath}
   The expression inside the absolute value is $3$ when
   $a_{00}=a_{01}=a_{10}=1$ and $a_{11}=0$. It cannot attain the
   value $\pm2$, since if so $2$ would factor into a product of
   four integers whose pairwise differences are all even. Hence
   $\displaystyle\lam(\z{2}^2)=\frac14\log3$.
\end{example}

Further numerical work suggests the following.

\begin{conjecture}\label{conj:Z/2's}
   $\displaystyle\lam(\z{2}^n)=\frac1{2^n}\log(2^n-1)$ for all $n\ge2$.
\end{conjecture}

Note that this quantity is the upper bound for $\lam(F)$ in Lemma
\ref{lem:upper-bound} when $F=\z{2}^n$.

\section{Mixed groups}\label{sec:mixed}

Here we consider groups $G$ that have a nontrivial connected
component $G^0$ of the identity.

Let $\Gamma=\ghat$, and denote the torsion subgroup of $\Gamma$
by $\tor(\Gamma)$. The totally disconnected group $G/G^0$ then
has dual group $\tor(\Gamma)$.

\begin{theorem}\label{thm:disconnected}
   If $G/G^0$ is infinite, then $\lam(G)=0$.
\end{theorem}

\begin{proof}
   Let $\Delta$ be a finite subgroup of $\tor(\Gamma)$, and
   suppose that $|\Delta|\ge3$. Then by
   Lemma \ref{lem:quotient},
   \begin{equation}\label{eqn:disconnected}
      \lam(G)\le\lam(G/\Delta^\perp)=\lam(\widehat{\Delta})\le
      \frac1{|\Delta|}\log\bigl(|\Delta|-1\bigr).
   \end{equation}
   If $\tor(\Gamma)$ is infinite, we can then find arbitrarily
   large finite subgroups $\Delta$, so that the bound in
   \eqref{eqn:disconnected} can be made arbitrarily small.
\end{proof}

\begin{remark}\label{rem:bound}
   If $G^0$ is nontrivial, then any character of infinite order
   that is nontrivial on $G^0$ gives a quotient map of $G$ to
   $\TT$, and so $\lam(G)\le\lam(\TT)$ by Lemma
   \ref{lem:quotient}. If $\lam(\TT)$ were $0$, this would show
   that only finite groups have positive Lehmer constant.
\end{remark}

\begin{example}\label{exam:mixed}
   What is the value of
   $\lam(\TT\oplus\z{2})$? If $\lam(\TT)$ were $0$, the previous
   remark shows that $\lam(\TT\oplus\z{2})$ would also be $0$.

   On the other hand, if $\lam(\TT)>0$, it seems likely that
   $\lam(\TT\oplus\z{2})<\lam(\TT)$. Some evidence for this comes
   from the following example, kindly communicated to us by Peter
   Borwein.

   Let
   \begin{displaymath}
      f(x)=x^{12}-x^{11}+x^{10}-x^9-x^6-x^3+x^2-x+1
   \end{displaymath}
   and
   \begin{displaymath}
      g(x)=x^8-x^7+x^6-x^5+x^4.
   \end{displaymath}
   Then $f+g$ is cyclotomic, so that $\m(f+g)=0$, while
   $\m(f-g)\approx 0.30082$.
      
   Define $h$ on $\TT\oplus\z{2}$ by
   \begin{displaymath}
      h(s,j)=f(e^{2\pi i s})+(-1)^j g(e^{2\pi is}).
   \end{displaymath}
   Then
   \begin{displaymath}
      \begin{aligned}
         \m_{\TT\oplus\z{2}}(h) &=\frac12 \bigl[
         \m(f+g)+\m(f-g) \bigr] \\
         &\approx 0.15041 < 0.16235 = \m_{\TT}(f_L),
      \end{aligned}
   \end{displaymath}
   where $f_L$ is the Lehmer polynomial in
   \eqref{eqn:lehmer-polynomial}. Thus $\lam(\TT\oplus\z{2})$ is
   strictly less than the best current value for $\lam(\TT)$.
\end{example}

\end{document}